\documentclass[11pt]{article}

\def\mathbb{\Bbb}
\usepackage{amsfonts,latexsym,amstext}
\usepackage{amsmath}
\usepackage{amssymb}

\usepackage{comment}
\usepackage{amsthm}
\usepackage[latin1]{inputenc}
\usepackage{color}

\theoremstyle{plain}
\newtheorem{theorem}{Theorem}[section]
\newtheorem{lemma}[theorem]{Lemma}

\newtheorem{definition}{Definition}[section]

\theoremstyle{remark}
\newtheorem{remark}[theorem]{Remark}

\makeatletter
\@addtoreset{equation}{section}

\makeatother

\def\qed{{\hfill\hbox{\enspace${ \square}$}} \smallskip}
\def\sqr#1#2{{\vcenter{\vbox{\hrule height .#2pt \hbox{\vrule
 width .#2pt height#1pt \kern#1pt \vrule
width .#2pt} \hrule height .#2pt}}}}
\def\square{\mathchoice\sqr54\sqr54\sqr{4.1}3\sqr{3.5}3}

\def\ds{\begin{displaystyle}}
\def\eds{\end{displaystyle}}

\def\<{\langle }
\def\>{\rangle }

\def\R{\mathbb R}

\def\E{\mathbb E}
\def\P{\mathbb P}

\def\H{\mathbb H}

\newcommand{\sper}[1]{\mathbb{E} \left[ #1 \right]}                               

\DeclareMathAlphabet{\mathonebb}{U}{bbold}{m}{n}                           %

\title{Existence and uniqueness for backward stochastic differential equations driven by a random measure}

\author{Elena Bandini\thanks{Politecnico di Milano, Dipartimento di Matematica, via Bonardi 9, 20133 Milano, Italy; ENSTA ParisTech, Unit\'e de Math\'ematiques appliqu\'ees, 828, boulevard des Mar\'echaux, F-91120 Palaiseau, France; e-mail: elena.bandini@polimi.it.}}

\date{}

\begin{document}
\allowdisplaybreaks

\maketitle

\begin{abstract}
\noindent We study the following backward stochastic differential equation on finite time horizon driven by an integer-valued random measure $\mu$ on $\R_+\times E$, where $E$ is a Lusin space, with compensator $\nu(dt,dx)=dA_t\,\phi_t(dx)$:
\[
Y_t = \xi + \int_{(t,T]} f(s,Y_{s-},Z_s(\cdot))\, d A_s - \int_{(t,T]} \int_E Z_s(x) \, (\mu-\nu)(ds,dx),\qquad 0\leq t\leq T.
\]
The generator $f$ satisfies, as usual, a uniform Lipschitz condition with respect to its last two arguments. In the literature, the existence and uniqueness for the above equation in the present general setting has only been established when $A$ is continuous or deterministic. The general case, i.e. $A$ is a right-continuous nondecreasing predictable process, is addressed in this paper. These results are relevant, for example, in the study of control problems related to Piecewise Deterministic Markov Processes (PDMPs). Indeed, when $\mu$ is the jump measure of a PDMP, then $A$ is predictable (but not deterministic) and discontinuous, with jumps of size equal to 1.

\vspace{4mm}


\vspace{4mm}

\end{abstract}

\section{Introduction}

Backward Stochastic Differential Equations (BSDEs) have been deeply studied since the seminal paper \cite{PardouxPeng}. In \cite{PardouxPeng}, as well as in many subsequent papers, the driving term was a Brownian motion. BSDEs with a discontinuous driving term have also been studied, see, among others, \cite{BP}, \cite{TangLi}, \cite{BBP}, \cite{ElKarouiHuang}, \cite{Xia}, \cite{Becherer}, \cite{CarboneFerrSantacroce}, \cite{CohenElliott}, \cite{JeanblancManiaSantacroceSchweizer}, \cite{ConfFuhrmanJacod}.

In all the papers cited above, and more generally in the literature on BSDEs, the generator (or driver) of the backward stochastic differential equation, usually denoted by $f$, is integrated with respect to a measure $dA$, where $A$ is a nondecreasing continuous (or deterministic and right-continuous as in \cite{CohenElliott}) process. The general case, i.e. $A$ is a right-continuous nondecreasing predictable process, is addressed in this paper. It is worth mentioning that Section 4.3 in \cite{ConfFuhrmanJacod} provides a \emph{counter-example} to existence for such general backward stochastic differential equations. For this reason, the existence and uniqueness result (Theorem \ref{T:MainThm}) is not a trivial extension of known results. Indeed, in Theorem \ref{T:MainThm} we have to impose an additional technical assumption, which is violated by the counter-example presented in \cite{ConfFuhrmanJacod} (see Remark \ref{R:Counter-example}(ii)). This latter assumption reads as follows: there exists $\varepsilon\in(0,1)$ such that (notice that $\Delta A_t \leq 1$)
\begin{equation}
\label{MainHyp_Intro}
2\,L_y^2\,|\Delta A_t|^2 \ \leq \ 1 - \varepsilon, \qquad \P\text{-a.s.},\,\forall\,t\in[0,T],
\end{equation}
where $L_y$ is the Lipschitz constant of $f$ with respect to $y$. As mentioned earlier, in \cite{CohenElliott} the authors study a class of BSDEs with a generator $f$ integrated with respect to a \emph{deterministic} (rather than predictable) right-continuous nondecreasing process $A$, even if this class is driven by a countable sequence of square-integrable martingales, rather than just a random measure. They provide an existence and uniqueness result for this class of BSDEs, see Theorem 6.1 in \cite{CohenElliott}, where the same condition \eqref{MainHyp_Intro} is imposed (see Remark \ref{R:Counter-example}(i)). However, the proof of Theorem 6.1 in \cite{CohenElliott} relies heavily on the assumption that $A$ is deterministic, and it can not be extended to the case where $A$ is predictable, which therefore requires a completely different proof.

As an application of the results presented in this paper, suppose that $\mu$ is the jump measure of a Piecewise Deterministic Markov Process (PDMP). Then, $A$ is predictable (not deterministic) and discontinuous, with jumps of size equal to 1. In this case condition \eqref{MainHyp_Intro} can be written as
\begin{equation}
\label{MainHypPDMP}
L_y \ < \ \frac{1}{\sqrt{2}}.
\end{equation}
This is the only additional condition required by Theorem \ref{T:MainThm}. In particular, Theorem \ref{T:MainThm} does not impose any condition on $L_z$, i.e. on the Lipschitz constant of $f$ with respect to its last argument. This is particularly important in the study of control problems related to PDMPs by means of BSDE methods. In this case $L_y=0$ and condition \eqref{MainHypPDMP} is automatically satisfied.

The paper is organized as follows: in Section \ref{Sec_preliminaries} we introduce the random measure $\mu$ and we fix the notation. In Section \ref{S:BSDE} we provide the definition of solution to the backward stochastic differential equation and we solve it in the case where $f=f(t,\omega)$ is independent of $y$ and $z$ (Lemma \ref{L:FirstCase}). Finally, in Section \ref{S:Main} we prove the main result (Theorem \ref{T:MainThm}) of this paper, i.e. the existence and uniqueness for our backward stochastic differential equation.

\section{Preliminaries}\label{Sec_preliminaries}

Consider a finite time horizon $T\in(0,\infty)$, a Lusin space $(E,\mathcal E)$, and a filtered probability space $(\Omega,\mathcal F,(\mathcal F_t)_{t\geq0},\mathbb P)$, with $(\mathcal F_t)_{t\geq0}$ right-continuous. We denote by $\mathcal P$ the predictable $\sigma$-field on $\Omega\times[0,T]$. In the sequel, given a measurable space $(G,\mathcal G)$, we say that a function on the product space $\Omega\times[0,T]\times G$ is predictable if it is $\mathcal P\otimes\mathcal G$-measurable.

Let $\mu$ be an integer-valued random measure on $\R_+\times E$. In the sequel we use a martingale representation theorem for the random measure $\mu$, namely Theorem 5.4 in \cite{J}. For this reason, we suppose that $(\mathcal F_t)_{t\geq0}$ is the natural filtration of $\mu$, i.e. the smallest right-continuous filtration in which $\mu$ is optional. We also assume that $\mu$ is a discrete random measure, i.e. the sections of the set $D=\{(\omega,t)\colon\mu(\omega,\{t\}\times E)=1\}$ are finite on every finite interval. However, the results of this paper (in particular, Theorem \ref{T:MainThm}) are still valid for more general random measure $\mu$ for which a martingale representation theorem holds (see Remark \ref{R:General_mu} for more details).

We denote by $\nu$ the $(\mathcal F_t)_{t\geq0}$-compensator of $\mu$. Then, $\nu$ can be disintegrated as follows
\[
\nu(\omega,dt,dx) \ = \ dA_t(\omega)\,\phi_{\omega,t}(dx),
\]
where $A$ is
a right-continuous nondecreasing predictable process such that $A_0=0$, and $\phi$ is a transition probability from $(\Omega\times[0,T],\mathcal P)$ into $(E,\mathcal E)$. We suppose, without loss of generality, that $\nu$ satisfies $\nu(\{t\}\times dx)\leq1$ identically, so that $\Delta A_t\leq1$. We define $A^c$ as $A_t^c=A_t-\sum_{0<s\leq t}\Delta A_s$, $\nu^c(dt,dx)=1_{J^c\times E}\,\nu(dt,dx)$, $\nu^d(dt,dx)=\nu(dt,dx)-\nu^c(dt,dx)=1_{J\times E}\,\nu(dt,dx)$, where $J=\{(\omega,t)\colon\nu(\omega,\{t\}\times dx)>0\}$.

\vspace{1mm}

We denote by $\mathcal B(E)$ the set of all Borel measurable functions on $E$. Given a measurable function $Z\colon\Omega\times[0,T]\times E\rightarrow\R$, we write $Z_{\omega,t}(x)=Z(\omega,t,x)$, so that $Z_{\omega,t}$, often abbreviated as $Z_t$ or $Z_t(\cdot)$, is an element of $\mathcal B(E)$. For any $\beta\geq0$ we also denote by $\mathcal E^\beta$ the Dol\'eans-Dade exponential of the process $\beta A$, which is given by
\begin{equation}
\label{ExpIncresing}
\mathcal{E}_t^{\beta} \ = \ e^{\beta\,A_t}\prod_{0 < s \leq t}(1 + \beta\,\Delta A_s)\,e^{-\beta\,\Delta\,A_s}.
\end{equation}

\section{The backward stochastic differential equation}
\label{S:BSDE}

The backward stochastic differential equation driven by the random measure $\mu$ is characterized by a triple $(\beta,\xi,f)$, where $\beta>0$ is a positive real number, and:
\begin{itemize}
\item $\xi\colon\Omega\rightarrow\R$, the \emph{terminal condition}, is an $\mathcal F_T$-measurable random variable satisfying $\E[\mathcal E_T^\beta |\xi|^2]$ $<\infty$;
\item $f\colon\Omega\times[0,T]\times\R\times\mathcal B(E)\rightarrow\R$, the \emph{generator}, is such that:
\begin{itemize}
\item[(i)] for any $y\in\R$ and $Z\colon\Omega\times[0,T]\times E\rightarrow\R$ predictable $\Longrightarrow$ $f(\omega,t,y,Z_{\omega,t}(\cdot))$ predictable;
\item[(ii)] for some nonnegative constants $L_y,L_z$, we have
\begin{align}
\label{Lipschitz_f}
&|f(\omega,t,y',\zeta') - f(\omega,t,y,\zeta)| \leq L_y|y'-y| \notag \\
&+ L_z\bigg(\int_E \bigg|\zeta'(x)-\zeta(x) - \Delta A_t(\omega)\int_E\big(\zeta'(z) - \zeta(z)\big)\phi_{\omega,t}(dz)\bigg|^2\,\phi_{\omega,t}(dx) \notag \\
&+ \Delta A_t(\omega)\big(1-\Delta A_t(\omega)\big)\bigg|\int_E (\zeta'(x)-\zeta(x))\,\phi_{\omega,t}(dx)\bigg|^2\bigg)^{1/2},
\end{align}
for all $(\omega,t)\in\Omega\times[0,T]$, $y,y'\in\R$, $\zeta,\zeta'\in L^2(E,\mathcal E,\phi_{\omega,t}(dx))$;
\item[(iii)] $\E[(1 + \sum_{0<t\leq T} |\Delta A_t|^2)\int_0^{T} \mathcal{E}_t^{\beta} |f(t,0,0)|^2 \,d A_t]<\infty$.
\end{itemize}
\end{itemize}
Given $(\beta,\xi,f)$, the backward stochastic differential equation takes the following form
\begin{equation}
\label{BSDE}
Y_t = \xi + \int_{(t,T]} f(s,Y_{s-},Z_s(\cdot))\, d A_s - \int_{(t,T]} \int_E Z_s(x) \, (\mu-\nu)(ds,dx),\qquad 0\leq t\leq T.
\end{equation}
\begin{definition}
\label{D:Space}
For every $\beta\geq0$, we define $\mathbb H_\beta^2(0,T)$ as the set of pairs $(Y,Z)$ such that:\begin{itemize}
\item $Y\colon\Omega\times[0,T]\rightarrow\R$ is an adapted c\`adl\`ag process satisfying
\begin{equation}
\label{Ybeta}
\|Y\|_{\H_{\beta,Y}^2(0,T)}^2 \ := \ \E\bigg[\int_{(0,T]} \mathcal E_t^\beta |Y_{t-}|^2\,dA_t\bigg] < \infty;
\end{equation}
\item $Z\colon\Omega\times[0,T]\times E\rightarrow\R$ is a predictable process satisfying
\begin{align}
\label{Zbeta}
\|Z\|_{\H_{\beta,Z}^2(0,T)}^2 \ &:= \ \E\bigg[\int_{(0,T]}\mathcal E_t^\beta\int_E \big|Z_t(x) - \hat Z_t\big|^2\,\nu(dt,dx) \notag \\
&\quad \ + \sum_{0 < t \leq T} \mathcal E_t^\beta\big|\hat Z_t\big|^2\big(1 - \Delta A_t\big)\bigg] \ < \ \infty,
\end{align}
where
\[
\hat Z_t \ = \ \int_E  Z_t(x)\,\nu(\{t\}\times dx), \qquad 0\leq t\leq T.
\]
\end{itemize}
For every $(Y,Z)\in\H_\beta^2(0,T)$, we denote
\[
\|(Y,Z)\|_{\mathbb H_\beta^2(0,T)}^2 := \|Y\|_{\H_{\beta,Y}^2(0,T)}^2 + \|Z\|_{\H_{\beta,Z}^2(0,T)}^2.
\]
\end{definition}

\begin{remark}
\label{R:Notation}
(i) Notice that the space $\H_\beta^2(0,T)$, endowed with the topology induced by $\|\cdot\|_{\H_\beta^2(0,T)}$, is an Hilbert space, provided we identify pairs of processes $(Y,Z),(Y',Z')$ satisfying $\|(Y-Y',Z-Z')\|_{\H_\beta^2(0,T)}=0$.

\vspace{1mm}

\noindent(ii) Suppose that there exists $\gamma\in(0,1]$ such that $\Delta A_t\leq1-\gamma$, for all $t\in[0,T]$, $\P$-a.s.. Then $Z$ belongs to $\H_{\beta,Z}^2(0,T)$ if and only if $\sqrt{\mathcal E^\beta} Z$ is in $L^2(\Omega\times[0,T]\times E,\mathcal P\otimes\mathcal E,\P\otimes\nu(dt,dx))$, i.e.
\[
\E\bigg[\int_{(0,T]}\mathcal E_t^\beta\int_E \big|Z_t(x)\big|^2\,\nu(dt,dx)\bigg] \ < \ \infty.
\]
%
%
\qed
\end{remark}

\begin{definition}
\label{D:Solution}
A solution to equation \eqref{BSDE} with data $(\beta,\xi,f)$ is a pair $(Y,Z)\in\H_\beta^2(0,T)$ satisfying equation \eqref{BSDE}. We say that equation \eqref{BSDE} admits a unique solution if, given two solutions $(Y,Z),(Y',Z')\in\H_\beta^2(0,T)$, we have $(Y,Z)=(Y',Z')$ in $\H_\beta^2(0,T)$.
\end{definition}

\begin{remark}
\label{R:Z_Mart}
Notice that, given a solution $(Y,Z)$ to equation \eqref{BSDE} with data $(\beta,\xi,f)$, we have (recalling that $\beta\geq0$, so that $\mathcal E_t^\beta\geq1$)
\[
\E\bigg[\int_{(0,T]}\int_E \big|Z_t(x) - \hat Z_t\big|^2\,\nu(dt,dx) + \sum_{0 < t \leq T}\big|\hat Z_t\big|^2\big(1 - \Delta A_t\big)\bigg] = \|Z\|_{\H_{0,Z}^2(0,T)}^2 \leq \|Z\|_{\H_{\beta,Z}^2(0,T)}^2 < \infty.
\]
This implies that the process $(Z_t1_{[0,T]}(t))_{t\geq0}$ belongs to $\mathcal G^2(\mu)$, see (3.62) and Proposition 3.71-(a) in \cite{J79}. In particular, the stochastic integral $\int_{(t,T]} \int_E Z_s(x) \, (\mu-\nu)(ds,dx)$ in \eqref{BSDE} is well-defined, and the process $M_t:=\int_{(0,t]}\int_E Z_s(x)(\mu-\nu)(ds,dx)$, $t\in[0,T]$, is a square integrable martingale (see Proposition 3.66 in \cite{J79}).
\qed
\end{remark}

\begin{lemma}\label{L:FirstCase}
Consider a triple $(\beta,\xi,f)$ and suppose that $f=f(\omega,t)$ does not depend on $(y,\zeta)$. Then, there exists a unique solution $(Y,Z)\in\H_\beta^2(0,T)$ to equation \eqref{BSDE} with data $(\beta,\xi,f)$. Moreover, the following identity holds:
\begin{eqnarray}\label{IdentityBSDE}
&&\E\big[\mathcal{E}_{t}^{\beta}\,|Y_t|^2\big] + \beta \,\E\bigg[\int_{(t,T]}\mathcal{E}_s^{\beta}\,(1+\beta\Delta A_s)^{-1}\,|Y_{s-}|^2\,d A_s\bigg] \nonumber\\
&& + \; \E\bigg[\int_{(t,T]}\mathcal E_s^\beta\int_E \big|Z_s(x) - \hat Z_s\big|^2\,\nu(ds,dx) + \sum_{t < s \leq T} \mathcal E_s^\beta\big|\hat Z_s\big|^2\big(1 - \Delta A_s\big)\bigg] \nonumber\\
&& = \; \E\big[\mathcal{E}_{T}^{\beta}\,|\xi|^2\big] + 2\,\E\bigg[\int_{(t,T]}\mathcal{E}_s^{\beta}\,Y_{s-} \, f_s\,d A_s\bigg] - \E\bigg[\sum_{t<s\leq T}\,\mathcal{E}_s^{\beta}\,|f_s|^2 \,|\Delta A_s|^2\bigg],
\end{eqnarray}
for all $t\in[0,T]$.
\end{lemma}
\proof
\emph{Uniqueness.} It is enough to prove that equation \eqref{BSDE} with data $(\beta,0,0)$ has the unique (in the sense of Definition \ref{D:Solution}) solution $(Y,Z)=(0,0)$. Let $(Y,Z)$ be a solution to equation \eqref{BSDE} with data $(\beta,0,0)$. Since the stochastic integral in \eqref{BSDE} is a square integrable martingale (see Remark \ref{R:Z_Mart}), taking the conditional expectation with respect to $\mathcal F_t$ we obtain, $\P$-a.s., $Y_t=0$, for all $t\in[0,T]$. This proves the claim for the component $Y$ and shows that the martingale $M_t:=\int_{(0,t]}\int_E Z_s(x)(\mu-\nu)(ds,dx)=0$, $\P$-a.s., for all $t\in[0,T]$. Therefore, the predictable bracket $\langle M,M\rangle_T=0$, $\P$-a.s., where we recall that (see Proposition 3.71-(a) in \cite{J79})
\[
\langle M,M\rangle_T \ = \ \int_{(0,T]}\int_E \big|Z_t(x) - \hat Z_t\big|^2\,\nu(dt,dx) + \sum_{0 < t \leq T} \big|\hat Z_t\big|^2\big(1 - \Delta A_t\big).
\]
This concludes the proof, since $\|Z\|_{\H_{\beta,Z}^2(0,T)}^2\leq\E\big[\mathcal E_T^\beta \langle M,M\rangle_T\big]=0$.

\vspace{1mm}

\noindent\emph{Identity \eqref{IdentityBSDE}.} Let $(Y,Z)$ be a solution to equation \eqref{BSDE} with data $(\beta,\xi,f)$. From It\^o's formula applied to
$\mathcal{E}_{s}^{\beta}\,|Y_s|^2$ it follows that (recall that $d\mathcal{E}_{s}^{\beta}= \beta\,\mathcal{E}_{s-}^{\beta}\,dA_s$)
\begin{eqnarray}\label{ItoFormula}
d(\mathcal{E}_{s}^{\beta}\,|Y_s|^2) &=&   \mathcal{E}_{s-}^{\beta}\,d|Y_{s}|^2 + |Y_{s-}|^2\,d\mathcal{E}_{s}^{\beta} + \Delta \mathcal{E}_{s}^{\beta}\,\Delta| Y_s|^2\notag\\
&=&   \mathcal{E}_{s-}^{\beta}\,d|Y_{s}|^2 + |Y_{s-}|^2\,d\mathcal{E}_{s}^{\beta} + (\mathcal{E}_{s} - \mathcal E_{s-}^{\beta})\,d| Y_s|^2\notag\\
&=&   \mathcal{E}_s^{\beta}\,d|Y_{s}|^2 + |Y_{s-}|^2\,d\mathcal{E}_{s}^{\beta} \notag\\
&=& 2\,\mathcal{E}_s^{\beta}\,Y_{s-}\,d Y_s + \mathcal{E}_s^{\beta}\,(\Delta Y_s)^2+\beta \, \mathcal{E}_{s-}^{\beta}\,|Y_{s-}|^2\,d A_s \notag \\
&=& 2\,\mathcal{E}_s^{\beta}\,Y_{s-}\,d Y_s + \mathcal{E}_s^{\beta}\,(\Delta Y_s)^2+\beta \, \mathcal E_s^\beta\,(1+\beta\Delta A_s)^{-1}\,|Y_{s-}|^2\,d A_s,
\end{eqnarray}
where the last equality follows from the identity $\mathcal E_{s-}^\beta=\mathcal E_s^\beta(1+\beta\Delta A_s)^{-1}$. Integrating \eqref{ItoFormula} on the interval $[t,\,T]$, we obtain
\begin{align}\label{expA_Ito}
\mathcal{E}_{t}^{\beta}\,|Y_t|^2 \ &= \
\mathcal{E}_{T}^{\beta}\,|\xi|^2 + 2 \int_{(t,T]}\mathcal{E}_s^{\beta}\,Y_{s-} \, f_s\,d A_s - 2\int_{(t,T]}\mathcal{E}_s^{\beta}\,Y_{s-}\int_{E}\,Z_{s}(x)\,(\mu - \nu)(ds,dx) \\
&\quad \ - \sum_{t<s\leq T}\,\mathcal{E}_s^{\beta}\,(\Delta Y_s)^2 - \beta \,\int_{(t,T]} \mathcal E_s^\beta\,(1+\beta\Delta A_s)^{-1}\,|Y_{s-}|^2\,d A_s. \notag
\end{align}
Now, notice that
\begin{equation}\label{Yjumps}
\Delta Y_s = \int_{E}\,Z_{s}(x)\,(\mu - \nu)(\{s\}\times dx) -f_s\,\Delta A_s.
\end{equation}
Thus
\begin{align}
|\Delta Y_s|^2 &= \bigg|\int_{E}\,Z_{s}(x)\,(\mu - \nu)(\{s\}\times dx)\bigg|^2 + |f_s|^2|\Delta A_s|^2 \notag \\
&\quad - 2 f_s\Delta A_s \int_{E}\,Z_{s}(x)\,(\mu - \nu)(\{s\}\times dx). \label{Yjumps_square}
\end{align}
Plugging \eqref{Yjumps_square} into \eqref{expA_Ito}, we find
\begin{eqnarray}\label{identity}
&&\mathcal{E}_{t}^{\beta}\,|Y_t|^2 + \beta \,\int_{(t,T]}\mathcal{E}_s^{\beta}\,(1+\beta\Delta A_s)^{-1}\,|Y_{s-}|^2\,d A_s + \sum_{t<s\leq T}\,\mathcal{E}_s^{\beta}\,\bigg|\int_{E}\,Z_{s}(x)\,(\mu - \nu)(\{s\}\times dx)\bigg|^2 \nonumber\\
&& = \; \mathcal{E}_{T}^{\beta}\,|\xi|^2 + 2 \int_{(t,T]}\mathcal{E}_s^{\beta}\,Y_{s-} \, f_s\,d A_s - 2\int_{(t,T]}\mathcal{E}_s^{\beta}\,Y_{s-}\int_{E}\,Z_{s}(x)\,(\mu - \nu)(ds,dx) \notag \\
&& - \; \sum_{t<s\leq T}\,\mathcal{E}_s^{\beta}\,|f_s|^2 \,|\Delta A_s|^2 + 2 \sum_{t<s\leq T}\,\mathcal{E}_s^{\beta}\,f_s\,\Delta A_s \int_{E}\,Z_{s}(x)\,(\mu - \nu)(\{s\}\times dx).
\end{eqnarray}
Notice that
\begin{align}
\label{NormZ}
&\E\bigg[\sum_{t<s\leq T}\,\mathcal{E}_s^{\beta}\,\bigg|\int_{E}\,Z_{s}(x)\,(\mu - \nu)(\{s\}\times dx)\bigg|^2\bigg] \notag \\
&= \E\bigg[\int_{(t,T]}\mathcal E_s^\beta\int_E \big|Z_s(x) - \hat Z_s\big|^2\,\nu(ds,dx) + \sum_{t < s \leq T} \mathcal E_s^\beta\big|\hat Z_s\big|^2\big(1 - \Delta A_s\big)\bigg].
\end{align}
We also observe that the two stochastic integrals
\begin{align*}
M^1_t &:=\int_{(0,t]}\mathcal{E}_s^{\beta}\,Y_{s-}\int_{E}Z_{s}(x)\,(\mu - \nu)(ds,dx) \\
M^2_t &:=\sum_{0<s\leq t}\mathcal{E}_s^{\beta}\,f_{s}\,\Delta A_s\int_{E}Z_{s}(x)\,(\mu - \nu)(\{s\}\times dx)
\end{align*}
are martingales. Therefore, taking the expectation in \eqref{identity} and using \eqref{NormZ}, we end up with \eqref{IdentityBSDE}. 

\vspace{1mm}

\noindent\emph{Existence.} Consider the martingale $\tilde M_t := \E[\xi + \int_{(0,T]} f_s\, d A_s|\mathcal{F}_t]$, $t\in[0,T]$. Let $M$ be a right-continuous modification of $\tilde M$. Then, by the martingale representation Theorem 5.4 in \cite{J} and Proposition 3.66 in \cite{J79} (noting that $M$ is a square integrable martingale), there exists a predictable process $Z\colon\Omega\times[0,T]\times E\rightarrow\R$ such that
\[
\E\bigg[\int_{(0,T]}\int_E \big|Z_t(x) - \hat Z_t\big|^2\,\nu(dt,dx) + \sum_{0 < t \leq T}\big|\hat Z_t\big|^2\big(1 - \Delta A_t\big)\bigg] < \infty
\]
and
\begin{equation}\label{rep1}
M_t = M_0 + \int_{(0,t]} \int_{E}Z_{s}(x)\, (\mu - \nu)(ds,dx), \qquad t \in [0,\,T].
\end{equation}
Set
\begin{equation}
\label{E:Ydef}
Y_t=M_t - \int_{(0,t]} f_s\, d A_s, \qquad t\in[0,T].
\end{equation}
Using the representation \eqref{rep1} of $M$, and noting that $Y_T =\xi$, we see that $Y$ satisfies \eqref{BSDE}.
When $\beta>0$, it remains to show that $Y$ satisfies \eqref{Ybeta} and $Z$ satisfies \eqref{Zbeta}. To this end,  let us define the increasing sequence of stopping times
\begin{eqnarray*}
	&&S_k = \inf \Big\{ t \in (0,\,T]: \int_{(0,t]}\mathcal{E}_{s}^{\beta}\,|Y_{s-}|^2\,d A_s  \\
	&&\qquad\qquad +\int_{(0,t]}\mathcal E_s^\beta\int_E \big|Z_s(x) - \hat Z_s\big|^2\,\nu(ds,dx) + \sum_{0 < s \leq t}\mathcal E_s^\beta\big|\hat Z_s\big|^2\big(1 - \Delta A_s\big)  > k \Big\}
\end{eqnarray*}
with the convention $\inf \emptyset = T$.
Computing the It\^o differential $d(\mathcal{E}_{s}^{\beta}\,|Y_s|^2)$ on the interval $[0,\,S_k]$ and proceeding as in the derivation of identity \eqref{IdentityBSDE}, we find
\begin{align}\label{expbeta_Ito1}
&\sper{\int_{(0,S_k]}\mathcal E_s^\beta\int_E \big|Z_s(x) - \hat Z_s\big|^2\,\nu(ds,dx) + \sum_{0 < s \leq S_k}\mathcal E_s^\beta\big|\hat Z_s\big|^2\big(1 - \Delta A_s\big)}\nonumber\\
	&+ \beta \,\sper{\int_{(0,S_k]}\mathcal{E}_{s}^{\beta}\,(1 + \beta \Delta A_s)^{-1}\,|Y_{s-}|^2\,d A_s} \leq \sper{\mathcal{E}_{S_k}^{\beta}|Y_{S_k}|^2}  + 2 \sper{\int_{(0,S_k]}\mathcal{E}_{s}^{\beta}\,Y_{s-} \, f_s\,d A_s}.
\end{align}
Let us now prove the following inequality (recall that we are assuming $\beta>0$)
\begin{equation}
\label{f}
\mathcal E_t^\beta \bigg(\int_{(t,T]} |f_s| \, dA_s\bigg)^2 \ \leq \ \bigg(\frac{1}{\beta} + \beta \sum_{t<s\leq T} |\Delta A_s|^2\bigg)\int_{(t,T]}\mathcal E_s^\beta\,|f_s|^2 \, dA_s.
\end{equation}
Set, for all $s\in[0,\,T]$,
\[
\bar A_s \ := \ \frac{\beta}{2}A_s^c + \sum_{0<r\leq s,\,\Delta A_r\neq0} \big(\sqrt{1 + \beta\Delta A_r} - 1\big), \qquad \underline A_s \ := \ - \frac{\beta}{2}A_s^c - \sum_{0<r\leq s,\,\Delta A_r\neq0} \frac{\sqrt{1 + \beta\Delta A_r} - 1}{\sqrt{1 + \beta\Delta A_r}}.
\]
Denote by $\bar{\mathcal E}$ (resp. $\underline{\mathcal E}$) the Dol\'eans-Dade exponential of the process $\bar A$ (resp. $\underline A$). Using Proposition 6.4 in \cite{J79} we see that
\begin{equation}
\label{Doleans-Dade}
1 \ = \ \underline{\mathcal E}_s\,\bar{\mathcal E}_s, \qquad (\bar{\mathcal E}_s)^2 \ = \ \mathcal E_s^\beta, \qquad \forall\,s\in[0,\,T].
\end{equation}
Then, we conclude that
\[
\mathcal E_t^\beta\bigg(\int_{(t,T]} |f_s| \, dA_s\bigg)^2 \ = \ \mathcal E_t^\beta\bigg(\int_{(t,T]} \underline{\mathcal E}_{s-}\,\bar{\mathcal E}_{s-}\,|f_s| \, dA_s\bigg)^2 \ \leq \ \bigg(\frac{1}{\beta} + \beta \sum_{t<s\leq T} |\Delta A_s|^2\bigg) \int_{(t,T]}\mathcal E_s^\beta\,|f_s|^2 \, dA_s,
\]
where we used the inequality $\mathcal E_{s-}^\beta\leq\mathcal E_s^\beta$ (which follows from \eqref{ExpIncresing}) and
\[
\mathcal E_t^\beta \int_{(t,T]}(\underline{\mathcal E}_{s-})^2\,dA_s \ = \ \mathcal E_t^\beta\frac{(\underline{\mathcal E}_t)^2 - (\underline{\mathcal E}_T)^2}{\beta} + \mathcal E_t^\beta \beta \sum_{t<s\leq T} (\underline{\mathcal E}_{s-})^2 \frac{|\Delta A_s|^2}{1 + \beta\,\Delta A_s} \ \leq \ \frac{1}{\beta} + \beta \sum_{t<s\leq T} |\Delta A_s|^2,
\]
where the last inequality follows from $\frac{1}{1+\beta\Delta A_s}\leq1$ and identities \eqref{Doleans-Dade}. Now, using \eqref{E:Ydef} and \eqref{f} we obtain
\begin{align}\label{estmate_expbeta_Y1}
	\mathcal{E}_t^{\beta}\,|Y_t|^2 \ &= \ \mathcal{E}_t^{\beta}\,\bigg|
	\mathbb E\bigg[\xi +\int_{(t,T]}\,f_s\,d A_s\Big|\mathcal{F}_t\bigg] \bigg|^2 \ \leq \ 2\,
	\mathbb E\big[\mathcal{E}_t^{\beta}\,|\xi|^2\big|\mathcal{F}_t\big]  +2\,\mathbb E\bigg[\mathcal{E}_t^{\beta}\,\bigg(\int_{(t,T]}\,|f_s|\,d A_s\bigg)^2\Big|\mathcal{F}_t\bigg] \nonumber\\
	&\leq \ 2\,
	\mathbb E\bigg[\mathcal{E}_T^{\beta}\,|\xi|^2 + \bigg(\frac{1}{\beta} + \beta \sum_{0<s\leq T} |\Delta A_s|^2\bigg)\int_{(0,T]}\mathcal{E}_s^{\beta}\,|f_s|^2\,d A_s\Big|\mathcal{F}_t\bigg].
\end{align}
Denote by $m_t$ a right-continuous modification of the right-hand side of \eqref{estmate_expbeta_Y1}. We see that $m=(m_t)_{t\in[0,T]}$ is a uniformly integrable martingale. In particular for every stopping time $S$ with values in $[0,\,T]$, we have, by Doob's optional stopping theorem,
\begin{equation}
\label{estmate_expbeta_Y1_bis}
	\sper{\mathcal{E}_{S}^{\beta}|Y_S|^2} \leq \sper{m_S} \leq \sper{m_T} < \infty.
\end{equation}
Notice that $(1+ \beta \Delta A_s)^{-1} \geq \frac{1}{1+\beta}$ $\P$-a.s. Using the inequality 2$ab \leq \gamma a^2 + \frac{1}{\gamma}b^2$ with $\gamma = \frac{\beta}{2(1+\beta)}$, 
and plugging \eqref{estmate_expbeta_Y1_bis} (with $S=S_k$) into \eqref{expbeta_Ito1}, 
we find the estimate
\begin{align*}
	&\frac{\beta}{2(1+\beta)}\sper{\int_{(0,S_k]} \mathcal{E}_{s}^{\beta} |Y_{s-}|^2\,d A_s} \\
	&+ \ \sper{\int_{(0,S_k]}\mathcal E_s^\beta\int_E \big|Z_s(x) - \hat Z_s\big|^2\,\nu(ds,dx) + \sum_{0 < s \leq S_k}\mathcal E_s^\beta\big|\hat Z_s\big|^2\big(1 - \Delta A_s\big)}\\
	&\leq \ 2 \,\mathbb E\big[\mathcal{E}_{T}^{\beta}\,|\xi|^2\big]+ 2\,
	\mathbb E\bigg[\bigg(\frac{1}{\beta} + \beta \sum_{0<s\leq T} |\Delta A_s|^2\bigg)\bigg(\int_{(0,T]} \mathcal{E}_{s}^{\beta} \,|f_s|^2\, d A_s\bigg)\bigg].
\end{align*}
From the above inequality we deduce that
\begin{eqnarray}\label{estimate_Sk1}
	&&\sper{\int_{(0,S_k]}\!\mathcal{E}_{s}^{\beta}\,|Y_{s-}|^2\,d A_s}
	+ \mathbb E\bigg[\int_{(0,S_k]}\!\mathcal E_s^\beta\!\int_E \big|Z_s(x) - \hat Z_s\big|^2\,\nu(ds,dx) + \!\sum_{0 < s \leq S_k}\!\mathcal E_s^\beta\big|\hat Z_s\big|^2\big(1 - \Delta A_s\big)\bigg]\nonumber\\
	&&\leq c(\beta)\,\left(\mathbb E\left[\mathcal{E}_{T}^{\beta}\,|\xi|^2\right]+\,\sper{\bigg(\frac{1}{\beta} + \beta \sum_{0<s\leq T} |\Delta A_s|^2\bigg)\int_{(0,T]} \mathcal{E}_{s}^{\beta} \,|f_s|^2\, d A_s}\right),
\end{eqnarray}
where $c(\beta) =2 + \frac{4(1 + \beta)}{\beta}$. Setting $S=\lim_{k}S_k$ we deduce
\begin{eqnarray*}
	&&\sper{\int_{(0,S]}\!\mathcal{E}_{s}^{\beta}\,|Y_{s-}|^2\,d A_s}+ 	+ \mathbb E\bigg[\int_{(0,S]}\!\mathcal E_s^\beta\int_E \big|Z_s(x) - \hat Z_s\big|^2\,\nu(ds,dx) + \sum_{0 < s \leq S}\!\mathcal E_s^\beta\big|\hat Z_s\big|^2\big(1 - \Delta A_s\big)\bigg]\\
	&&< \infty, \quad \P\text{-a.s.,}
\end{eqnarray*}
which implies $S=T$, $\P$-a.s., by the definition of $S_k$. Letting $k \rightarrow \infty$ in \eqref{estimate_Sk1}, we conclude that $Y$ satisfies \eqref{Ybeta} and $Z$ satisfies \eqref{Zbeta}, so that $(Y,Z)\in\H_\beta^2(0,T)$.
\endproof

\section{Main result}
\label{S:Main}

\begin{theorem}
\label{T:MainThm}
Suppose that there exists $\varepsilon\in(0,1)$ such that
\begin{equation}
\label{MainHyp}
2\,L_y^2\,|\Delta A_t|^2 \ \leq \ 1 - \varepsilon, \qquad \P\text{-a.s.},\,\forall\,t\in[0,T].
\end{equation}
Then there exists a unique solution $(Y,Z)\in\H_\beta^2(0,T)$ to equation \eqref{BSDE} with data $(\beta,\xi,f)$, for every $\beta$ satisfying
\begin{equation}
\label{IneqBeta}
\beta \ \geq \ \frac{\frac{L_y^2}{\hat L_{z,t}^2} + \frac{2\,\hat L_{z,t}^2}{1 - \delta + 2\,\hat L_{z,t}^2\,\Delta A_t}}{1 - \Delta A_t\Big(\frac{L_y^2}{\hat L_{z,t}^2} + \frac{2\,\hat L_{z,t}^2}{1 - \delta + 2\,\hat L_{z,t}^2\,\Delta A_t}\Big)}, \qquad \P\text{-a.s.},\,\forall\,t\in[0,T],
\end{equation}
for some $\delta\in(0,\varepsilon)$ and strictly positive predictable process $(\hat L_{z,t})_{t\in[0,T]}$ given by\begin{equation}
\label{L_z}
\hat L_{z,t}^2 \ = \ \max\bigg(L_z^2 + \delta,\frac{(1 - \delta)\,L_y}{\sqrt{2(1 - \delta)} - 2\,L_y\,\Delta A_t}\bigg).
\end{equation}
\end{theorem}
\begin{remark}
(i) Notice that when condition \eqref{MainHyp} holds the right-hand side of \eqref{IneqBeta} is a well-defined nonnegative real number, so that there always exists some $\beta\geq0$ which satisfies \eqref{IneqBeta}.

\vspace{1mm}

\noindent(ii) Observe that in Theorem \ref{T:MainThm} there is no condition on $L_z$, i.e. on the Lipschitz constant of $f$ with respect to its last argument.
\qed
\end{remark}
\proof[Proof of Theorem \ref{T:MainThm}]
The proof is based on a fixed point argument that we now describe.
Let us consider the function $\Phi: \H_\beta^2(0,T)\rightarrow\H_\beta^2(0,T)$, mapping
$(U,V)$ to $(Y,Z)$ as follows:
\begin{eqnarray}\label{fixed_point_eq}
Y_{t} &=& \xi + \int_{(t,T]} f(t,\,U_{s-},\,V_s)\, d A_s
- \int_{(t,T]} \int_{E} Z_{s}(x) \, (\mu - \nu)(ds,dx), \qquad 0\leq t\leq T.
\end{eqnarray}
By Lemma  \ref{L:FirstCase} there exists a unique  $(Y,Z)\in\H_\beta^2(0,T)$ satisfying \eqref{fixed_point_eq}, so that $\Phi$ is a well-defined map.
We then see that $(Y,Z)$ is a solution in $\H_\beta^2(0,T)$ to the BSDE \eqref{BSDE} with data $(\beta,\xi,f)$ if and only if it is a fixed point of $\Phi$.

Let us prove that $\Phi$ is a contraction when $\beta$ is large enough. Let $(U^i,V^i)\in\H_\beta^2(0,T)$, $i =1,2$, and set $(Y^i,Z^i) = \Phi(U^i,V^i)$. Denote $\bar{Y}= Y^1-Y^2$, $\bar{Z}= Z^1-Z^2$, $\bar{U}= U^1-U^2$, $\bar{V}= V^1-V^2$, $\bar{f}_{s} = f(s,\,U_{s-}^1,\,V_s^1) - f(s,\,U^2_{s-},\,V^2_s)$. Notice that
\begin{equation}
\label{BSDE_diff}
\bar Y_t \ = \ \int_{(t,T]} \bar f_s\, d A_s - \int_{(t,T]} \int_E \bar Z_s(x) \, (\mu-\nu)(ds,dx),\qquad 0\leq t\leq T.
\end{equation}
Then, identity \eqref{IdentityBSDE}, with $t=0$, becomes (noting that $\E[\mathcal{E}_{0}^{\beta}|\bar Y_0|^2]$ is nonnegative)
\begin{eqnarray}\label{IdentityBSDE2}
&&\beta \,\E\bigg[\int_{(0,T]}\mathcal{E}_s^{\beta}\,(1+\beta\Delta A_s)^{-1}\,|\bar Y_{s-}|^2\,d A_s\bigg] \nonumber\\
&& + \; \E\bigg[\int_{(0,T]}\mathcal E_s^\beta\int_E \big|\bar Z_s(x) - \hat{\bar Z}_s\big|^2\,\nu(ds,dx) + \sum_{0 < s \leq T} \mathcal E_s^\beta\big|\hat{\bar Z}_s\big|^2\big(1 - \Delta A_s\big)\bigg] \nonumber\\
&& \leq \; 2\,\E\bigg[\int_{(0,T]}\mathcal{E}_s^{\beta}\,\bar Y_{s-} \,\bar f_s\,d A_s\bigg] - \E\bigg[\sum_{0<s\leq T}\,\mathcal{E}_s^{\beta}\,|\bar f_s|^2 \,|\Delta A_s|^2\bigg].
\end{eqnarray}
From the standard inequality $2ab\leq\frac{1}{\alpha} a^2+\alpha b^2$, $\forall\,a,b\in\R$ and $\alpha>0$, we obtain, for any strictly positive predictable processes $(c_s)_{s\in[0,T]}$ and $(d_s)_{s\in[0,T]}$,
\begin{align*}
2\,\E\bigg[\int_{(0,T]}\mathcal{E}_s^{\beta}\,\bar Y_{s-} \,\bar f_s\,d A_s\bigg] &\leq \E\bigg[\int_{(0,T]}\frac{1}{c_s}\mathcal{E}_s^{\beta}\,|\bar Y_{s-}|^2\,d A_s^c\bigg] + \E\bigg[\sum_{0<s\leq T}\frac{1}{d_s}\,\mathcal{E}_s^{\beta}\,|\bar Y_{s-}|^2\,\Delta A_s\bigg] \\
&\quad \ + \E\bigg[\int_{(0,T]}c_s\,\mathcal{E}_s^{\beta}\,|\bar f_s|^2\,d A_s^c\bigg] + \E\bigg[\sum_{0<s\leq T}d_s\,\mathcal{E}_s^{\beta}\,|\bar f_s|^2\,\Delta A_s\bigg].
\end{align*}
Therefore \eqref{IdentityBSDE2} becomes
\begin{align}\label{IdentityBSDE3}
&\E\bigg[\int_{(0,T]} \Big(\beta - \frac{1}{c_s}\Big) \mathcal E_s^\beta |\bar Y_{s-}|^2\,dA_s^c\bigg] + \E\bigg[\sum_{0<s\leq T}\Big(\beta\,(1+\beta\Delta A_s)^{-1} - \frac{1}{d_s}\Big)\,\mathcal{E}_s^{\beta}\,|\bar Y_{s-}|^2\,\Delta A_s\bigg] \nonumber\\
& + \E\bigg[\int_{(0,T]}\mathcal E_s^\beta\int_E \big|\bar Z_s(x) - \hat{\bar Z}_s\big|^2\,\nu(ds,dx) + \sum_{0 < s \leq T} \mathcal E_s^\beta\big|\hat{\bar Z}_s\big|^2\big(1 - \Delta A_s\big)\bigg] \nonumber\\
& \leq \E\bigg[\int_{(0,T]}c_s\,\mathcal{E}_s^{\beta}\,|\bar f_s|^2\,d A_s^c\bigg] + \E\bigg[\sum_{0<s\leq T}\,\big(d_s - \Delta A_s\big)\,\mathcal{E}_s^{\beta}\,|\bar f_s|^2 \,\Delta A_s\bigg].
\end{align}
Now, by the Lipschitz property \eqref{Lipschitz_f} of $f$, we see that for any predictable process $(\hat L_{z,s})_{s\in[0,T]}$, satisfying $\hat L_{z,s}>L_z$, $\P$-a.s. for every $s\in[0,\,T]$, we have
\begin{equation}
\label{f_Lipschitz}
|\bar f_s|^2 \leq 2L_y^2|\bar U_{s-}|^2 + 2\hat L_{z,s}^2\bigg(\int_E \big|\bar V_s(x) - \hat{\bar V}_s\big|^2\,\phi_s(dx) + 1_{\{\Delta A_s\neq0\}}\frac{1 - \Delta A_s}{\Delta A_s}\big|\hat{\bar V}_s\big|^2\bigg),
\end{equation}
for all $s\in[0,\,T]$. For later use, fix $\delta\in(0,\varepsilon)$ and take $(\hat L_{z,s})_{s\in[0,T]}$ given by \eqref{L_z}. Notice that the two components inside the maximum in \eqref{L_z} are nonnegative (the first being always strictly positive, the second being zero if $L_y=0$) and uniformly bounded, as it follows from condition \eqref{MainHyp}. Plugging inequality \eqref{f_Lipschitz} into \eqref{IdentityBSDE3}, and using the following identity for $\bar Z$ (and the analogous one for $\bar V$)
\begin{align*}
&\E\bigg[\int_{(0,T]}\mathcal E_s^\beta\int_E \big|\bar Z_s(x) - \hat{\bar Z}_s\big|^2\,\nu(ds,dx) + \sum_{0 < s \leq T} \mathcal E_s^\beta\big|\hat{\bar Z}_s\big|^2\big(1 - \Delta A_s\big)\bigg] \\
&= \ \E\bigg[\int_{(0,T]} \mathcal E_s^\beta \int_E |\bar Z_s(x)|^2\,\nu^c(ds,dx)\bigg] + \E\bigg[\sum_{0 < s \leq T} \mathcal E_s^\beta \big(\widehat{|\bar Z_s|^2} - |\hat{\bar Z}_s|^2\big)\bigg],
\end{align*}
we obtain
\begin{align}\label{IdentityBSDE4}
&\E\bigg[\int_{(0,T]} \Big(\beta - \frac{1}{c_s}\Big)\mathcal E_s^\beta |\bar Y_{s-}|^2\,dA_s^c\bigg] + \E\bigg[\sum_{0<s\leq T}\Big(\beta\,(1+\beta\Delta A_s)^{-1} - \frac{1}{d_s}\Big)\,\mathcal{E}_s^{\beta}\,|\bar Y_{s-}|^2\,\Delta A_s\bigg] \nonumber\\
& + \E\bigg[\int_{(0,T]} \mathcal E_s^\beta \int_E |\bar Z_s(x)|^2\,\nu^c(ds,dx)\bigg] + \E\bigg[\sum_{0 < s \leq T} \mathcal E_s^\beta \big(\widehat{|\bar Z_s|^2} - |\hat{\bar Z}_s|^2\big)\bigg] \nonumber\\
& \leq \; 2\,L_y^2\,\E\bigg[\int_{(0,T]} c_s\,\mathcal E_s^\beta |\bar U_{s-}|^2\,dA_s^c\bigg] + 2\,\E\bigg[\int_{(0,T]} c_s\,\hat L_{z,s}^2\, \mathcal E_s^\beta \int_E |\bar V_s(x)|^2\,\nu^c(ds,dx)\bigg] \\
& + 2\,L_y^2\,\E\bigg[\sum_{0<s\leq T}\big(d_s - \Delta A_s\big)\,\mathcal{E}_s^{\beta}\,|\bar U_{s-}|^2 \,\Delta A_s\bigg] + 2\,\E\bigg[\sum_{0<s\leq T}\big(d_s - \Delta A_s\big)\,\hat L_{z,s}^2\,\mathcal{E}_s^{\beta}\,\big(\widehat{|\bar V_s|^2} - |\hat{\bar V}_s|^2\big) \bigg]. \notag
\end{align}
Set $b_s:=\min(\beta-\frac{1}{c_s},\beta(1+\beta\Delta A_s)^{-1}-\frac{1}{d_s})$ and $a_s:=2\hat L_{z,s}^2\max(c_s,d_s-\Delta A_s)$, $s\in[0,\,T]$. Then, inequality \eqref{IdentityBSDE4} can be rewritten as (recalling that $\hat L_{z,s}>0$)
\begin{align}\label{IdentityBSDE5}
&\E\bigg[\int_{(0,T]} b_s\,\mathcal E_s^\beta |\bar Y_{s-}|^2\,dA_s^c\bigg] + \E\bigg[\sum_{0<s\leq T}\,b_s\,\mathcal{E}_s^{\beta}\,|\bar Y_{s-}|^2\,\Delta A_s\bigg] \nonumber\\
& + \E\bigg[\int_{(0,T]} \mathcal E_s^\beta \int_E |\bar Z_s(x)|^2\,\nu^c(ds,dx)\bigg] + \E\bigg[\sum_{0 < s \leq T} \mathcal E_s^\beta \big(\widehat{|\bar Z_s|^2} - |\hat{\bar Z}_s|^2\big)\bigg] \nonumber\\
& \leq \ \E\bigg[ \int_{(0,T]} \frac{L_y^2}{\hat L_{z,s}^2}\,a_s\, \mathcal E_s^\beta |\bar U_{s-}|^2\,dA_s^c\bigg] + \E\bigg[\sum_{0<s\leq T}\frac{L_y^2}{\hat L_{z,s}^2}\,a_s\,\mathcal{E}_s^{\beta}\,|\bar U_{s-}|^2 \,\Delta A_s\bigg] \notag\\
& + \E\bigg[ \int_{(0,T]} a_s\,\mathcal E_s^\beta \int_E |\bar V_s(x)|^2\,\nu^c(ds,dx)\bigg] + \E\bigg[\sum_{0<s\leq T}\,a_s\,\mathcal{E}_s^{\beta}\,\big(\widehat{|\bar V_s|^2} - |\hat{\bar V}_s|^2\big) \bigg]. 
\end{align}
It follows from \eqref{IdentityBSDE5} that $\Phi$ is a contraction if:
\begin{enumerate}
\item[(i)] there exists $\alpha\in(0,1)$ such that $a_s\leq\alpha$, $\P$-a.s. for every $s\in[0,\,T]$;
\item[(ii)] $\frac{L_y^2}{\hat L_{z,s}^2}\leq b_s$, $\P$-a.s. for every $s\in[0,\,T]$.
\end{enumerate}
Let us prove that (i) and (ii) hold. Regarding (i), we have, for all $s\in[0,\,T]$,
\[
c_s \ \leq \ \frac{1 - \alpha}{2\,\hat L_{z,s}^2}, \qquad d_s \ \leq \ \frac{1 - \alpha}{2\,\hat L_{z,s}^2} + \Delta A_s.
\]
It is useful for condition (ii) to choose $\alpha=\delta$, where $\delta\in(0,\varepsilon)$ was fixed in the statement of the theorem, and $c_s,d_s$ given by
\begin{equation}
\label{d_s}
c_s \ = \ \frac{1 - \delta}{2\,\hat L_{z,s}^2}, \qquad d_s \ = \ \frac{1 - \delta}{2\,\hat L_{z,s}^2} + \Delta A_s,
\end{equation}
for all $s\in[0,\,T]$. Concerning (ii), we have, for all $s\in[0,\,T]$,
\[
\min\Big(\beta-\frac{1}{c_s},\beta(1+\beta\Delta A_s)^{-1}-\frac{1}{d_s}\Big) \ \geq \ \frac{L_y^2}{\hat L_{z,s}^2},
\]
which becomes
\begin{equation}
\label{beta}
\beta \ \geq \ \frac{L_y^2}{\hat L_{z,s}^2} + \frac{1}{c_s}, \qquad \beta \ \geq \ \frac{\frac{L_y^2}{\hat L_{z,s}^2} + \frac{1}{d_s}}{1 - \Delta A_s\Big(\frac{L_y^2}{\hat L_{z,s}^2} + \frac{1}{d_s}\Big)},
\end{equation}
where for the last inequality we need to impose the additional condition
\[
1 - \Delta A_s\bigg(\frac{L_y^2}{\hat L_{z,s}^2} + \frac{1}{d_s}\bigg) \ > \ 0.
\]
This latter inequality can be rewritten as
\begin{equation}
\label{LyDeltaA}
L_y^2\,\Delta A_s \ < \ \hat L_{z,s}^2\bigg(1 - \frac{\Delta A_s}{d_s}\bigg) \ = \ \frac{(1 - \delta)\,\hat L_{z,s}^2}{1 - \delta + 2\,\hat L_{z,s}^2\,\Delta A_s},
\end{equation}
where the last equality follows from the definition of $d_s$ in \eqref{d_s}. From \eqref{L_z}, and since in particular
\[
\hat L_{z,s}^2 \ \geq \ \frac{(1 - \delta)\,L_y}{\sqrt{2(1 - \delta)} - 2\,L_y\,\Delta A_s} \ > \ \frac{(1 - \delta)\,L_y^2\,\Delta A_s}{1 - \delta - 2\,L_y^2\,|\Delta A_s|^2}, \qquad \P\text{-a.s.},\,\forall\,s\in[0,\,T],
\]
it follows that inequality \eqref{LyDeltaA} holds. Finally, concerning \eqref{beta}, we begin noting that
\[
\frac{L_y^2}{\hat L_{z,s}^2} + \frac{1}{c_s} < \frac{\frac{L_y^2}{\hat L_{z,s}^2} + \frac{1}{d_s}}{1 - \Delta A_s\Big(\frac{L_y^2}{\hat L_{z,s}^2} + \frac{1}{d_s}\Big)},
\]
as it can be shown using \eqref{d_s}. Now, let us denote
\[
\frac{\frac{L_y^2}{\hat L_{z,s}^2} + \frac{1}{d_s}}{1 - \Delta A_s\Big(\frac{L_y^2}{\hat L_{z,s}^2} + \frac{1}{d_s}\Big)} \ = \ H_s(\hat L_{z,s}^2),
\]
where, for every $s\in[0,\,T]$,
\[
H_s(\ell) \ = \ \frac{h_s(\ell)}{1 - \Delta A_s\,h_s(\ell)}, \qquad h_s(\ell) \ = \ \frac{L_y^2}{\ell} + \frac{2\,\ell}{1 - \delta + 2\,\ell\,\Delta A_s}, \qquad \ell>0.
\]
Notice that $H_s$ attains its minimum at
$\ell_s^*=\frac{(1 - \delta)\,L_y}{\sqrt{2(1 - \delta)} - 2\,L_y\,\Delta A_s}$. This explains the expression of the second component inside the maximum in \eqref{L_z}. In conclusion, given $(\hat L_{z,s})_{s\in[0,\,T]}$ as in \eqref{L_z} we obtain a lower bound for $\beta$ from the second inequality in \eqref{beta}, which corresponds to \eqref{IneqBeta}.
\endproof

\begin{remark}
\label{R:Counter-example}
(i) In \cite{CohenElliott} the authors study a class of BSDEs driven by a countable sequence of square-integrable martingales, with a generator $f$ integrated with respect to a right-continuous nondecreasing process $A$ as in \eqref{BSDE}. Similarly to our setting, $A$ is not necessarily continuous, however in \cite{CohenElliott} it is supposed to be \emph{deterministic} (instead of predictable). Theorem 6.1 in \cite{CohenElliott} provides an existence and uniqueness result for the class of BSDEs studied in \cite{CohenElliott} under the following assumption ($2\,L_{y,t}^2$ corresponds to $c_t$ and $\Delta A_t$ corresponds to $\Delta\mu_t$ in the notation of \cite{CohenElliott}):
\begin{equation}
\label{CohenElliott}
2\,L_{y,t}^2\,|\Delta A_t|^2 \ < \ 1, \qquad \forall\,t\in[0,\,T],
\end{equation}
where $L_{y,t}$ is a measurable deterministic function uniformly bounded such that \eqref{Lipschitz_f} holds with $L_{y,t}$ in place of $L_y$. As showed at the beginning of the proof of Theorem 6.1 in \cite{CohenElliott}, if \eqref{CohenElliott} holds (and $A$ is as in \cite{CohenElliott}), then there exists $\varepsilon\in(0,1)$ such that
\begin{equation}
\label{MainHyp_bis}
2\,L_{y,t}^2\,|\Delta A_t|^2 \ \leq \ 1 - \varepsilon, \qquad \forall\,t\in[0,\,T].
\end{equation}
This proves that when condition \eqref{CohenElliott} holds then \eqref{MainHyp_bis} is also valid, since in our setting we can take $L_{y,t}\equiv L_y$.

\vspace{1mm}

\noindent(ii) Section 4.3 in \cite{ConfFuhrmanJacod} provides a counter-example to existence for BSDE \eqref{BSDE} when $A$ is discontinuous, as it can be the case in our setting; the rest of the paper \cite{ConfFuhrmanJacod} studies BSDE \eqref{BSDE} with $A$ continuous. Let us check that the counter-example proposed in \cite{ConfFuhrmanJacod} does not satisfy condition \eqref{MainHyp}. In \cite{ConfFuhrmanJacod} the process $A$ is a pure jump process with a single jump of size $p\in(0,1)$ at a deterministic time $t\in(0,T]$. The Lipschitz constant of $f$ with respect to $y$ is $L_y=\frac{1}{p}$. Then
\[
2\,L_y^2\,|\Delta A_t|^2 \ = \ 2
\]
if $t$ is the jump time of $A$, so that condition \eqref{MainHyp} is violated.\qed
\end{remark}

\begin{remark}
\label{R:General_mu}
Suppose that $\mu$ is an integer-valued random measure on $\R_+\times E$ not necessarily discrete. Then $\nu$ can still be disintegrated as follows
\[
\nu(\omega,dt,dx) \ = \ dA_t(\omega)\,\phi_{\omega,t}(dx),
\]
where $A$ is a right-continuous nondecreasing predictable process such that $A_0=0$, but $\phi$ is in general only a transition measure (instead of transition probability) from $(\Omega\times[0,T],\mathcal P)$ into $(E,\mathcal E)$. Notice that when $\mu$ is discrete one can choose $\phi$ to be a transition probability, therefore $\phi(E)=1$ and $\nu(\{t\}\times E)=\Delta A_t$ (a property used in the previous sections). When $\mu$ is not discrete, let us suppose that $\nu^d$ can be disintegrated as follows 
\begin{equation}
\label{nu_decomposition}
\nu^d(\omega,dt,dx) \ = \ \Delta A_t(\omega)\,\phi_{\omega,t}^d(dx), \qquad \phi_{\omega,t}^d(E) \ = \ 1,
\end{equation}
where $\phi^d$ is a transition \emph{probability} from $(\Omega\times[0,T],\mathcal P)$ into $(E,\mathcal E)$. In particular $\nu^d(\{t\}\times E)=\Delta A_t$. Then, when \eqref{nu_decomposition} and a martingale representation theorem for $\mu$ hold, all the results of this paper are still valid and can be proved proceeding along the same lines. As an example, \eqref{nu_decomposition} holds when $\mu$ is the jump measure of a L\'evy process, indeed in this case $\Delta A_t$ is identically zero.
\qed
\end{remark}

\paragraph{Acknowledgements.} The author would like to thank Prof. Jean Jacod for his helpful discussions and valuable suggestions to improve this paper.

\end{document}